# Upper bounds for the travel time on traffic systems

Nadir Farhi[*], Habib Haj-Salem, and Jean-Patrick Lebacque

*Université Paris-Est, IFSTTAR/COSYS/GRETTIA, F-77447 Champs-sur Marne Cedex France.*

**Abstract.** A key measure of performance and comfort in a road traffic network is the travel time that the users of the network experience to complete their journeys. Travel times on road traffic networks are stochastic, highly variable, and dependent on several parameters. It is, therefore, necessary to have good indicators and measures of their variations. In this article, we extend a recent approach for the derivation of deterministic bounds on the travel time in a road traffic network (Farhi, Haj-Salem and Lebacque 2013). The approach consists in using an algebraic formulation of the cell-transmission traffic model on a ring road, where the car-dynamics is seen as a linear min-plus system. The impulse response of the system is derived analytically, and is interpreted as what is called a service curve in the network calculus theory (where the road is seen as a server). The basic results of the latter theory are then used to derive an upper bound for the travel time through the ring road.

We consider in this article open systems rather than closed ones. We define a set of elementary traffic systems and an operator for the concatenation of such systems. We show that the traffic system of any road itinerary can be built by concatenating a number of elementary traffic systems. The concatenation of systems consists in giving a service guarantee of the resulting system in function of service guarantees of the composed systems. We illustrate this approach with a numerical example, where we compute an upper bound for the travel time on a given route in a urban network.

## Introduction

We are interested in this article in the derivation of upper bounds for the travel time through a given route in a road network. The approach we adopt here, is algebraic. It has been introduced in(Farhi, Haj-Salem and Lebacque 2013). In the latter reference, a such bound is derived from an algebraic formulation of a first order traffic model on a single-lane ring road. Moreover, the bound on the travel time is given in function of the average car-density in the ring road. The derivation performed in (Farhi, Haj-Salem and Lebacque 2013) is based on basic results of the network calculus theory (Chang 2000) (Le Boudec and Thiran 2001) (Cruz 1991). More precisely, it is first shown that the traffic dynamics derived from the model can be written linearly in a specific algebraic structure (the min-plus algebra of functions; (Baccelli et al. 1992)(Chang 2000)). Then, from that formulation, the impulsion response of the min-plus linear system is interpreted as what is called a *service curve* of the road (seen as a server) in the network calculus theory. Finally, for a given traffic inflow to the road, and by using the service curve derived before, a basic result of the network calculus gives an upper bound of the travel time through the ring road. See (Farhi, Haj-Salem and Lebacque 2013) for more details. The application of the Network Calculus theory to road traffic control has been recently treated in (Varaiya 2013). The approach we adopt here bases on the same theory, but is different from that of (Varaiya 2013).

---

[*] Corresponding author. e-mail: `nadir.farhi@ifsttar.fr`

In this article we follow the same ideas as in (Farhi, Haj-Salem and Lebacque 2013), but we adopt a system theory approach. It consists in defining a number of elementary traffic systems, for which service curves can be calculated, with adequate algebraic operators for the connection of those systems. Therefore, one can build large scale road networks by connecting predefined elementary traffic systems. Moreover, we give a result that tells how to derive a service curve of a system resulting from the connection of two systems with known service curves. Thus, a service curve of a road network can be derived, and from that, upper bounds for the travel time through the network routes.

**Notations**

| | |
|---|---|
| $F$ | set of non decreasing non negative time functions (cumulated flows). |
| $\oplus$ | addition in $F$. such that $(f \oplus g)(t) = \min(f(t), g(t))$. |
| $*$ | product (min. convolution) in $F$. such that $(fg)(t) = (f*g)(t) = \min_{0 \leq s \leq t}(f(s) + g(t-s))$. |
| $(F, \oplus, *)$ | idempotent semi-ring. |
| $\varepsilon$ | zero element for $(F, \oplus, *)$, such that $\varepsilon(t) = +\infty, \forall t \geq 0$. |
| $e$ | unity element for $(F, \oplus, *)$, such that $e(0) = 0, e(t) = +\infty, \forall t > 0$. |
| $\gamma^p$ | gain signal in $F$. such that $\gamma^p(0) = p$ and $\gamma^p(t) = +\infty$ for $t > 0$. |
| $\delta^T$ | shift signal in $F$. such that $\delta^T(t) = 0$ for $0 \leq t \leq T$ and $\delta^T(t) = +\infty$ for $t > T$. |
| $\oslash$ | substraction (min. de-convolution) in $F$. such that $(f \oslash g)(t) = \max_{s \geq 0}(f(t+s) - g(s))$. |
| $f^k$ | power convolution in $F$. such that, for $f \in F, f^0 = e$, and $f^k = f^{k-1} * f$, for $k \geq 1$. |
| $f^*$ | sub-additive closure in $F$. such that for $f \in F, f^* = \oplus_{k \geq 0} f^k$. |
| $M$ | set of $n \times m$ matrices with entries in $F$. |
| $\oplus$ | addition in $M$. such that $(A \oplus B)_{ij}(t) = (A_{ij} \oplus B_{ij})(t) = \min(A_{ij}(t), B_{ij}(t))$. |
| $*$ | product in $M$. such that $(AB)_{ij}(t) = (A*B)_{ij}(t) = \min_{1 \leq k \leq n}(A_{ik} * B_{kj})(t)$. |
| $(M_{n \times n}, \oplus, *)$ | idempotent semi-ring. |
| $\varepsilon$ | zero element for $(M_{n \times n}, \oplus, *)$, such that $\varepsilon_{ij} = \varepsilon, \forall i, j$. |
| $e$ | unity element for $(M_{n \times n}, \oplus, *)$, such that $e_{ii} = e, \forall i$, and $e_{ij} = \varepsilon, \forall i \neq j$. |
| $A^k$ | power operation in $M_{n \times n}$. such that $A^0 = e$ and $A^k = A^{k-1} * A$, for $k \geq 1$. |
| $A^*$ | sub-additive closure in $M_{n \times n}$; such that $A^* = \oplus_{k \geq 0} A^k$. |

# 1. Review and complements in network calculus

We give in this section a short review of the main definitions and results that we use in this article. The main variables we use here are the cumulated traffic flows, which we denote with capital letters function of time. Any traffic model is then seen as a system with input signals (car inflows) and output signals (car outflows). The *network calculus* theory associates an *arrival* and a *service curves* to a system, and derives from those curves performance bounds like upper bounds of the delay of passing through the system. An arrival curve consists in upper-bounding the arrival inflows to the system, while a service curve consists in lower-bounding the guaranteed service and then the departure outflows from the system. In the following, we review these two notions of arrival and service curve in the one-dimensional case (for systems with one arrival inflow, and one departure outflow. As mentioned in the table of notations above, we consider the set $F$ of non-decreasing and non-negative time functions, in which we consider the two operations: addition $\oplus$ (element-wise minimum) and the product $*$ (minimum convolution); see the notations above. Let us now, consider a system (seen as a server) with an arrival cumulated flow $U \in F$, and a departure cumulated flow $Y \in F$.

**Definition 1**. A curve $\alpha$ is an arrival curve for $U$, if $U \leq \alpha * U$.

**Definition 2.** A curve $\beta$ is a service curve for the server, if $Y \geq \beta * U$.

For the derivation of performance bounds from arrival and service curves in the one-dimension case, see (Cruz 1991) (Chang 2000), (Le Boudec and Thiran 2001).

We are concerned here by the multi-dimensional case, where multiple inflows arrive to and departure from the system. The particular signals $\varepsilon, e, \gamma^p$ and $\delta^T$, defined in the table of notations above, will be used here. We consider the following definitions and results.

**Definition 3**. (Arrival matrix). For a given $n \times 1$ vector $U$ of cumulated arrival flows $U_i, i = 1, \cdots, n$, a $n \times n$ matrix $\alpha$ is said to be a $T$-arrival matrix for $U$ if

$$\forall i,j = 1,2,\cdots,n, \qquad U_i \leq \delta^{-T_{ij}} \alpha_{ij} \otimes U_j$$

That is to say that

$$\forall i,j = 1,2,\cdots,n, \forall s,t \in \mathbb{N}, \qquad U_i(t) - U_j(s) \leq \alpha_{ij}(T_{ij} + t - s). \blacksquare$$

Let us notice that for $i = j$, the curve $\alpha_{ii}$ is simply a one-dimension arrival curve for $U_i$. However, for $i \neq j$, the difference with respect to the case $i = j$ is that, it is possible to have, $U_i(t) - U_j(s) > 0$, even for $t < s$. Indeed, if we assume that $U_i(t) - U_j(s) \leq 0, \forall t < s$, then we get $U_i(s) - U_j(s) \leq 0, \forall s \geq 0$, and similarly $U_j(s) - U_i(s) \leq 0, \forall s \geq 0$. Therefore, $U_i(s) - U_j(s) = 0, \forall s \geq 0$. It is trivial that such an assumption is very restrictive. Therefore, if we like to upper bound $U_i(t) - U_j(s)$ for all $s, t \geq 0$, then we need to work with negative times for the arrival curve. In order to continue working with non-negative times, we time-shift here the curve with negative times to zero.

A simple way to obtain such T-arrival matrices, is, first to determine the matrix $T$ (of non negative entries). For a given couple $(i, j)$, $T_{ij}$ is determined as follows.

$$T_{ij} = \text{Min}\{\tau \geq 0, U_i(t+\tau) - U_j(t) \geq 0, \forall t \geq 0\}.$$

Then, $\alpha_{ij}$ is determined using Definition 3:

$$\alpha_{ij} \geq \delta^{T_{ij}} U_i \oslash U_j$$

It is easy to check that for $i = j$, we have $T_{ii} = 0$, and then $\alpha_{ii}$ is a one-dimensional arrival curve for $U_i$.

Let us notice that Definition 3 is different form Definition 4.2.1 given in (Chang 2000). Definition 3 is illustrated in the numerical example of the last section.

**Definition 4.** (Service matrix). For a given server with input vector $U$ and output vector $Y$, a $n \times n$ matrix $\beta$ is said to be a service matrix for the server, if $Y \geq \beta * U$.

Definition 5. (Virtual delay). For a given server with input vector $U$ and output vector $Y$, the virtual delay of the last quantity arrived at time $t$ from the $i^{\text{th}}$ input to depart from the $i^{\text{th}}$ output, denoted $d_i(t)$ is defined:

$$d_i(t) := \inf\{d \geq 0, Y_i(t+d) \geq U_i(t)\}.$$

Theorem 1. For a given server with input vector $U$ and output vector $Y$, if $\alpha$ is a T-arrival matrix for $U$, and $\beta$ is a service matrix for the server, then

$$\forall i = 1,2,\cdots,n, \forall t \in \mathbb{N}, d_i(t) \leq \text{Inf}\{d \geq 0, \alpha_{ij}(T_{ij} + s) \leq \beta_{ij}(s + d),$$
$$-T_{ij} \leq s \leq t, j = 1,2,\cdots,n\}.$$

and then the virtual delays $d_i, i = 1,\cdots,n$ are bounded as follows.

$$\forall i = 1,2,\cdots,n, \forall t \in \mathbb{N}, \quad d_i(t) \leq \max_{1 \leq j \leq n}\{\sup_{s \geq -T_{ij}}\{inf\{d \geq 0, \alpha_{ij}(T_{ij} + s) \leq \beta_{ij}(s + d)\}\}\}.$$

Proof. It is a trivial adaptation of the proof of Theorem 4.3.6 in (Chang 2000). ∎

In order to build a whole traffic network, we base in elementary traffic systems which will be used as unit systems in the composition. We consider here two elementary systems (an uncontrolled and a controlled road sections). The composition we use here are inspired from (Farhi 2008.); see also (Farhi 2009), (Farhi, Goursat and Quadrat 2001), (Farhi, Goursat and Quadrat 2005), (Farhi, Goursat and Quadrat 2007), (Farhi, Goursat and Quadrat 2011) and (Farhi 2012).

## 2. The road section model

We consider a road section system $i$, as illustrated in Figure 1. Cars arrive from the left side of the road section, pass through it, and departs from the right side of it. The inputs $U_{fw}$ and $U_{bw}$ represent respectively the traffic demand from the upstream section $i - 1$ to the section $i$, and the traffic supply of the downstream section $i + 1$ to the section $i$. The outputs $Y_{fw}$ and $Y_{bw}$ represent respectively the traffic demand of the section $i$ to the downstream section $i + 1$, and the traffic supply of the section $i$ to the upstream section $i - 1$.

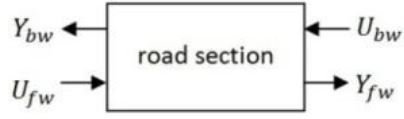

Figure 1. A road section.

Let us clarify the notations $U_{fw}, Y_{fw}, U_{bw}$ and $Y_{bw}$:
- $U_{fw}(t)$ (resp. $Y_{fw}$): cumulated forward inflow (resp. outflow) of cars from time zero to time $t$.
- $U_{bw}(t)$ (resp. $Y_{bw}$): cumulated backward supply of section $i+1$ (resp. $i$) from time zero to time $t$.

Let us define the cumulated flow $Z_{fw}$ as follows.

$$Z_{fw}(t) := \left(Y_{fw}(t) - n\right)^+ := \max(Y_{fw}(t) - n, 0),$$

where $n$ denotes the number of cars in the considered section at time zero. Therefore, under the condition $U_{fw}(0) = Y_{fw}(0)$, and under FIFO condition, the cumulated flow $Z_{fw}$ coincides with the outflow corresponding to the inflow $U_{fw}$; that is to say that at any time $t$, $Z_{fw}(t)$ gives the cumulated number of cars that leaved the road section up to time $t$, among from the cumulated number $U_{fw}(t)$ of cars arrived to the road section up to time $t$.

Similarly, we define the cumulated flow $Z_{bw}$ as follows.

$$Z_{bw}(t) = (Y_{bw}(t) - \bar{n})^+ := \max(Y_{bw}(t) - n, 0),$$

where $\bar{n} = n_{max} - n$, with $n_{max}$ denoting the maximum number of cars that the road section can contain. Again similarly, under the condition $U_{bw}(0) = Y_{bw}(0)$, and under FIFO condition, the cumulated flow $Z_{bw}$ coincides with the outflow corresponding to the inflow $U_{bw}$.

In order to simplify the presentation of the ideas, and without loss of generality, we assume that all the cumulated flows are initialized to zero.

$$U_{fw}(0) = Y_{fw}(0) = U_{bw}(0) = Y_{bw}(0) = 0, \tag{1}$$

Let us now write the traffic dynamics on the road section. As in (Farhi, Haj-Salem and Lebacque 2013), we base on the cell-transmission model (Daganzo 1994) with a trapezoidal fundamental diagram; see also (Lebacque 1996). We obtain the following dynamics, where we introduce an intermediate variable $Q$, which is simply the cumulated forward outflow $Y_{fw}$.

$$\begin{aligned}
Q(t) &= \min\left\{Q\left(t - \frac{\Delta x}{v}\right) + q_{max}\frac{\Delta x}{v}, U_{fw}\left(t - \frac{\Delta x}{v}\right) + n, U_{bw}(t)\right\} \\
Y_{fw}(t) &= Q(t) \\
Y_{bw}(t) &= Q\left(t - \frac{\Delta x}{w}\right) + \bar{n}
\end{aligned} \tag{2}$$

By using the min-plus algebra notations (see (Farhi, Haj-Salem and Lebacque 2013), (Baccelli et al. 1992)), we get:

$$Q = \gamma^{q_{max}\,\Delta x/v}\delta^{\Delta x/v}Q \oplus \gamma^n \delta^{\Delta x/v} U_{fw} \oplus U_{bw}$$
$$Y_{fw} = Q \oplus e$$
$$Y_{bw} = \gamma^{\bar{n}}\delta^{\Delta x/w} Q \oplus e$$

where we added (min-plus addition) the unity vector $e$ to $Y_{fw}$ and $Y_{bw}$ in order to satisfy condition (1). Then, by denoting $U = (U_{fw} \quad U_{bw})$, and $Y = (Y_{fw} \quad Y_{bw})$, we can write

$$Q = A * Q \oplus B * U$$
$$Y = C * Q \oplus e \qquad (3)$$

with $A = \gamma^{q_{max}\,\Delta x/v}\delta^{\Delta x/v}$, $B = (\gamma^n \delta^{\Delta x/v} \quad e)$, and $C = \begin{pmatrix} e \\ \gamma^{\bar{n}}\delta^{\Delta x/w} \end{pmatrix}$. Let us notice that in (3), $e$ denotes also the vector $(e \quad e)$.

Therefore, the traffic dynamics on a road section may be represented with three matrices: a signal $A$, a line matrix $B$ of two signals, and a column matrix $C$ of two signals. In a such configuration, $B$ represents the traffic demand of the section, $C$ gives the traffic supply that the section offers for an eventual upstream section, while $A$ gives the outflow limit $q_{max}$ imposed by the trapezoidal fundamental diagram.

From (3), (see (Baccelli et al. 1992) or (Chang 2000)), $Y \geq (C * A^* * B) * U \oplus e$. Then, since $U(0) = 0$, we have $e \geq U$. Hence

$$Y \geq (C * A^* * B) * U \oplus U = (e \oplus C * A^* * B) * U.$$

Then from the expressions of $Z_{fw}$ and $Z_{bw}$, we obtain

$$Z \geq H * (e \oplus C * A^* * B) * U. \qquad (4)$$

with

$$H = \begin{pmatrix} \gamma^{-n} & e \\ e & \gamma^{-\bar{n}} \end{pmatrix}.$$

One can easily check that (see (Farhi, Haj-Salem and Lebacque 2013)).

$$C * A^* * B = \left(\gamma^{q_{max}\,\Delta x/v}\delta^{\Delta x/v}\right)^* * \begin{pmatrix} \gamma^n \delta^{\Delta x/v} & e \\ \gamma^{n+\bar{n}}\delta^{\Delta x/v + \Delta x/w} & \gamma^{\bar{n}}\delta^{\Delta x/w} \end{pmatrix}.$$

**Theorem 2.** The matrix $H * (e \oplus C * A^* * B)$ is a service matrix for the road section, seen as a server with two inputs and two outputs.

**Corollary 1.** A service matrix $\beta$ for the road section seen as a server, is given as follows.

$$\beta = \begin{pmatrix} q_{max}\left(t - \dfrac{\Delta x}{v}\right)^+ & q_{max}\left(t - \dfrac{\Delta x}{v} - \dfrac{n}{q_{max}}\right)^+ \\ q_{max}\left(t - \dfrac{\Delta x}{v} - \dfrac{\Delta x}{v} + \dfrac{n}{q_{max}}\right)^+ & q_{max}\left(t - \dfrac{\Delta x}{w}\right)^+ \end{pmatrix}$$

**Proof.** It consists in proving that $H * (e \oplus C * A^* * B) \geq \beta$. ∎

## 3. The controlled road section model

We consider here a road section controlled with a traffic light. We denote by $c$ the cycle time of the traffic light, and by $G$ and $R$ the green and red times, with $c = G + R$. We can easily check that the traffic dynamics in the control road section is the same as (2), except the first equation that changes to

$$Q(t) = \min\left\{Q\left(t - \frac{\Delta x}{v}\right) + \left(\frac{G}{c}\right)q_{max}\frac{\Delta x}{v}, U_{fw}\left(t - \frac{\Delta x}{v} - R\right) + n, U_{bw}(t)\right\} \quad (5)$$

The first term in (5) replaces $Q\left(t - \frac{\Delta x}{v}\right) + q_{max}\frac{\Delta x}{v}$ in (2) since $G/c < 1$.

The term $U_{fw}\left(t - \frac{\Delta x}{v}\right) + n$ in (2) is time-shifted by $R$ in (5). Indeed cars arriving to the light may have an additional delay upper-bounded by $R$.

The dynamics (5) tells that the inflow to the traffic light passes through the light with a maximum time delay of $R$ time units, under the supply constraint downstream of the light, and with a maximum flow of $G/C$.

**Theorem 3**. The matrix $H * (e \oplus C * A'^* * B')$ is a service matrix for the road section, seen as a server with two inputs and two outputs, with $A' = \gamma^{(G/c)q_{max}\Delta x/v}\delta^{\Delta x/v}$, $B' = \begin{pmatrix}\gamma^n\delta^{\Delta x/v+R} & e\end{pmatrix}$.

**Proof.** Directly from Theorem 2, with the modifications made in (5). ∎

## 4. Concatenation of traffic systems

The composition of traffic systems is done in two dimensions, since each system has two inputs and two outputs. The connection is not in series, in the sense that connecting two systems does not mean connect the outputs of one system to the inputs of the other. Indeed, the connection is made here in the two directions, by connecting one output of system 1 to one input of system 2, and one output of system 2 to one input of system 1. In Figure 2, we illustrate the connection of two elementary systems (road sections).

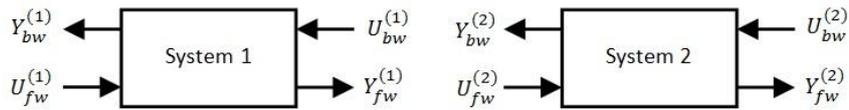

Figure 2. Concatenation of two min-plus linear traffic systems.

Let us consider two min-plus linear traffic systems 1 and 2, with service matrices $\beta^{(1)}$ and $\beta^{(2)}$. We then have:

$$\begin{pmatrix}Y_{fw}^{(i)}\\Y_{bw}^{(i)}\end{pmatrix} = \begin{pmatrix}(\beta^{(i)})_{11} & (\beta^{(i)})_{12}\\(\beta^{(i)})_{21} & (\beta^{(i)})_{22}\end{pmatrix}\begin{pmatrix}U_{fw}^{(i)}\\U_{bw}^{(i)}\end{pmatrix}, \quad i = 1,2.$$

The following result is on the composition such two systems.

**Theorem 4.** A service matrix $\beta$ for the whole system is given by:

$$\beta_{11} = \beta_{11}^{(2)}\beta_{11}^{(1)} \oplus \beta_{11}^{(2)}\beta_{12}^{(1)}\left(\beta_{21}^{(2)}\beta_{12}^{(1)}\right)^*\beta_{21}^{(2)}\beta_{11}^{(1)}$$
$$\beta_{12} = \beta_{11}^{(2)}\beta_{12}^{(1)}\left(\beta_{21}^{(2)}\beta_{12}^{(1)}\right)^*\beta_{22}^{(2)} \oplus \beta_{12}^{(2)}$$
$$\beta_{21} = \beta_{21}^{(1)} \oplus \beta_{22}^{(1)}\left(\beta_{21}^{(2)}\beta_{12}^{(1)}\right)^*\beta_{21}^{(2)}\beta_{11}^{(1)}$$
$$\beta_{22} = \beta_{22}^{(1)}\left(\beta_{21}^{(2)}\beta_{12}^{(1)}\right)^*\beta_{22}^{(2)}$$

such that

$$\begin{pmatrix} Y_{fw}^{(2)} \\ Y_{bw}^{(1)} \end{pmatrix} = \begin{pmatrix} \beta_{11} & \beta_{12} \\ \beta_{21} & \beta_{22} \end{pmatrix} * \begin{pmatrix} U_{fw}^{(1)} \\ U_{bw}^{(2)} \end{pmatrix}.$$

**Proof.** See Appendix A. ∎

## 5. Roads and itineraries

In order to build a road of $m$ sections, we need to compose $m$ elementary traffic systems of road sections. The service matrix of each road section can be obtained by Theorem 2, giving fundamental traffic diagrams on each section. Then the service matrix of the whole road is obtained by the composition of the road section systems and by applying Theorem 4. A controlled road of $m$ sections is obtained similarly by composing $m - 1$ uncontrolled road sections with one controlled road section.

A route (or an itinerary) in a controlled road network is build by composing a number of controlled roads. In Figure 5, we illustrate the composition of controlled roads to obtain a traffic system associated to a whole road network. The procedure of computing a service matrix for a traffic flow passing respectively through roads R1, R2, R3, and R4 is the following.

- Determine service matrices for all uncontrolled sections of the itinerary, by Theorem 2.
- Determine service matrices for all the controlled sections of the itinerary, by Theorem 3.
- Determine a service matrix for the itinerary by connecting the systems R1, R2, R3, R4, by Theorem 4.

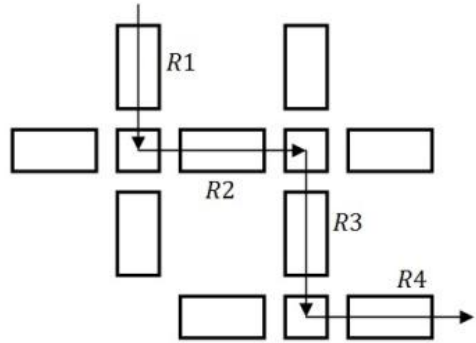

Figure 3. An itinerary of four roads and three intersections.

Once a service matrix is determined for a given traffic system, and having an arrival matrix expressing the traffic demand arriving to the system and the traffic supply that backwards to it, it suffices to apply Theorem 1 to obtain upper bounds for the travel time for any input - output couple of the traffic system.

## 6. A numerical example

We illustrate the results presented in this article with a numerical example. We consider the itinerary of Figure 3, for which we calculate an upper bound for the travel time from the entry of road 1 to the exit from road 4 passing through the roads R1, R2, R3 and R4. We consider some common parameters for all the road sections: $v = 15$ meter/sec., $w = 7$ m/s, $\rho_j = 1/10$ veh/m. Other parameters for the road sections are given in Table 1.

|  | R1 | R2 | R3 | R4 |
|---|---|---|---|---|
| Length $\Delta x$ (meter) | 150 | 150 | 100 | 100 |
| Maximum flow $q_{max}$ (veh/sec) | 0.32 | 0.35 | 0.4 | 0.38 |
| Initial density of cars $n$ (veh/meter) | 5/150 | 10/150 | 3/100 | 7/100 |
| Cycle time (sec.) | 60 | 90 | 80 | - |
| Green time (sec.) | 30 | 50 | 45 | - |

Figure 4. Parameters of the road sections R1, R2, R3 and R4 of Figure 3.

The results of this example are illustrated in Figure 4. The input signals $U_{fw}$ arriving to road 1, and $U_{bw}$ backing from road 4 are taken such that the arrival flows do not exceed the service offered by the whole route. The arrival curves of the arrival matrix $\alpha$ are computed by Definition 3. First the shift times $T_{12} = 60\ s.$, and $T_{21} = 8\ s.$ are computed. Then the curves are deduced by Definition 3. The service curves are computed following the steps cited above. An upper bound for the travel time through the route is then calculated according to Theorem 1. We are concerned here by the delay $d_1$ corresponding the forward travel time (the delay $d_2$ corresponds to the backward travel time of the backward waves). We obtained for this example the following result.

$$d_1 = \max(d_{11}, d_{12}) = \max(205, 241) = 241 \text{ seconds.}$$

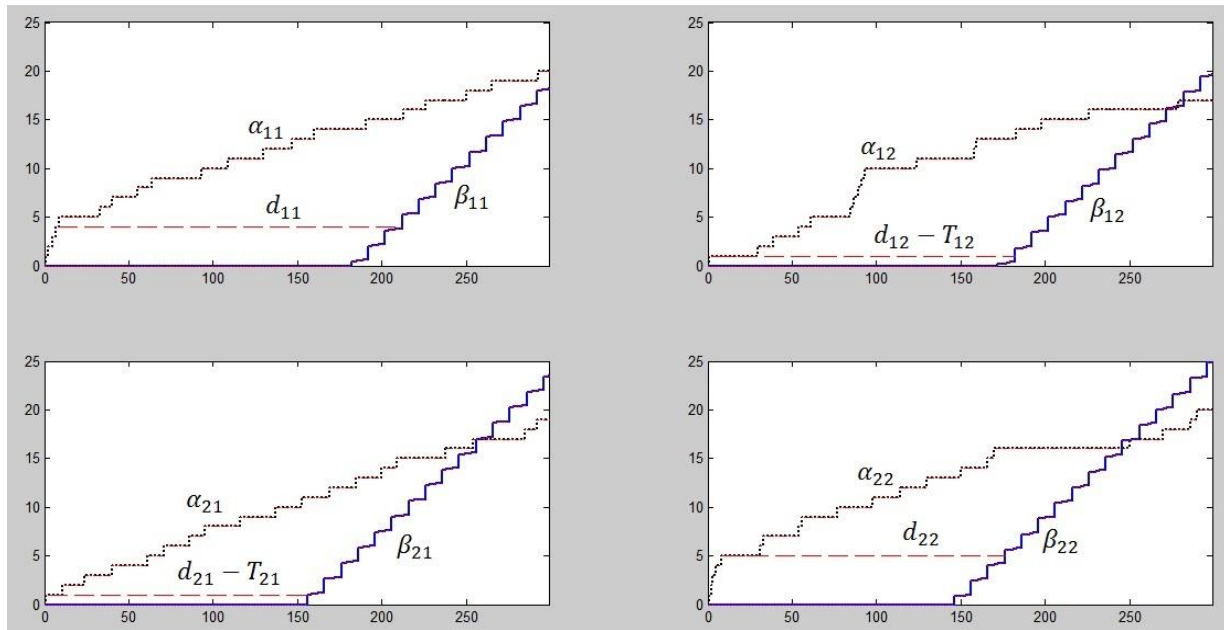

Figure 5. Arrival curves of the arrival matrix, service curves of the service matrix, and the time delays.

## References

BACCELLI, Francois, et al. (1992). *Synchronization and Linearity*. Wiley.
CHANG, Cheng-Shang (2000). *Performance guarantees in communication networks*. Springer.
CRUZ, R. (1991). A calculus of network delay. Part I: Network elements in isolation. Part II: Network analysis. **1** (137 pp. 114-141),.
DAGANZO, C. F. (1994). *The cell transmission model: A dynamic representation of highway traffic consistent with the hydrodynamic theory*. Transportation Research Part B: Methodological, 28(4), 269-287.
FARHI, Nadir (2009). Modeling and control of elementary 2D-Traffic systems using Petri nets and minplus algebra. in Proceedings of the IEEE Conference on Decision and Control.
FARHI, Nadir (2008.). *Modélisation minplus et commande du trafic de villes régulières.*. PhD thesis, University of Paris 1 Panthéon - Sorbonne. ed.,
FARHI, Nadir (2012). Piecewise linear car-following modeling. *Transportation Research Part C.*, **25**, 100-112.
FARHI, Nadir, GOURSAT, Maurice and QUADRAT, Jean.-Pierre. (2005). Derivation of the fundamental traffic diagram for two circular roads and a crossing using minplus algebra and Petri net modeling.. In: *IEEE Conference on Decision and Control.*,.
FARHI, Nadir, GOURSAT, Maurice and QUADRAT, Jena-Pierre. (2007). Fundamental Traffic Diagrams of Elementary Road Networks. In: *European Control Conference.*,.
FARHI, Nadir, GOURSAT, Maurice and QUADRAT, Jean-Pierre. (2001). Piecewise linear concave dynamical systems appearing in the microscopic traffic modeling. *Linear Algebra and its Applications.*, **435** (7), 1711-1735.
FARHI, Nadir, GOURSAT, Maurice and QUADRAT, Jean-Pierre. (2011). The traffic phases of road networks. *Transportation Research Part C.*, **19** (1), 85-102.
FARHI, Nadir, HAJ-SALEM, Habib and LEBACQUE, Jean-Patrick. (2013). Algebraic approach for performance bound calculus in transportation networks (Road Network calculus).. *Transportation Research Record.*,.
LE BOUDEC, J.-Y. and THIRAN, P. (2001). *Network calculus - A theory for deterministic queuing systems for the internet.* Springer-Verlag..
LEBACQUE, J-P. (1996). *The Godunov scheme and what it means for the first order traffic flow models*. In: Lesort, J.-B. (Ed.), Proceedings of the 13th ISTTT, 647-678.
VARAIYA, Pravin (2013). *The Max-Pressure controller for arbitrary networks for signalised intersections.* Advances in Dynamic Network Modeling in Complex Transportation Systems, Springer.

**Proof of Theorem 4.**

From the dynamics of system 2, we have $Y_{bw}^{(2)} = \beta_{21}^{(2)} U_{fw}^{(2)} \oplus \beta_{22}^{(2)} U_{bw}^{(2)}$.
Then by replacing $U_{fw}^{(2)}$ by $Y_{fw}^{(1)}$, we get

$$\begin{aligned} Y_{bw}^{(2)} &= \beta_{21}^{(2)} \left[ \beta_{11}^{(1)} U_{fw}^{(1)} \oplus \beta_{12}^{(1)} U_{bw}^{(1)} \right] \oplus \beta_{22}^{(2)} U_{bw}^{(2)} \\ &= \beta_{21}^{(2)} \beta_{11}^{(1)} U_{fw}^{(1)} \oplus \beta_{21}^{(2)} \beta_{12}^{(1)} U_{bw}^{(1)} \oplus \beta_{22}^{(2)} U_{bw}^{(2)} \end{aligned}$$

Now, we replace $U_{bw}^{(1)}$ by $Y_{bw}^{(2)}$. We obtain

$$Y_{bw}^{(2)} = \beta_{21}^{(2)} \beta_{11}^{(1)} U_{fw}^{(1)} \oplus \beta_{21}^{(2)} \beta_{12}^{(1)} Y_{bw}^{(2)} \oplus \beta_{22}^{(2)} U_{bw}^{(2)},$$

for which the solution, in $Y_{bw}^{(2)}$ is given as follows.

$$Y_{bw}^{(2)} = \left( \beta_{21}^{(2)} \beta_{12}^{(1)} \right)^* \left( \beta_{21}^{(2)} \beta_{11}^{(1)} U_{fw}^{(1)} \oplus \beta_{22}^{(2)} U_{bw}^{(2)} \right). \tag{6}$$

From the dynamics of system 2, we also have

$$Y_{fw}^{(2)} = \beta_{11}^{(2)} U_{fw}^{(2)} \oplus \beta_{12}^{(2)} U_{bw}^{(2)}.$$

Then by replacing $U_{fw}^{(2)}$ by $Y_{fw}^{(1)}$, we get

$$\begin{aligned} Y_{fw}^{(2)} &= \beta_{11}^{(2)} \left( \beta_{11}^{(1)} U_{fw}^{(1)} \oplus \beta_{12}^{(1)} U_{bw}^{(1)} \right) \oplus \beta_{12}^{(2)} U_{bw}^{(2)} \\ &= \beta_{11}^{(2)} \beta_{11}^{(1)} U_{fw}^{(1)} \oplus \beta_{11}^{(2)} \beta_{12}^{(1)} U_{bw}^{(1)} \oplus \beta_{12}^{(2)} U_{bw}^{(2)} \end{aligned}$$

We then replace $U_{bw}^{(1)}$ by the expression of $Y_{bw}^{(2)}$ in (6). We get
$$Y_{fw}^{(2)} = \beta_{11}^{(2)} \beta_{11}^{(1)} U_{fw}^{(1)} \oplus \beta_{11}^{(2)} \beta_{12}^{(1)} \left( \beta_{21}^{(2)} \beta_{12}^{(1)} \right)^* \left( \beta_{21}^{(2)} \beta_{11}^{(1)} U_{fw}^{(1)} \oplus \beta_{22}^{(2)} U_{bw}^{(2)} \right) \oplus \beta_{12}^{(2)} U_{bw}^{(2)}.$$

That is

$$\begin{aligned} Y_{fw}^{(2)} = &\left[ \beta_{11}^{(2)} \beta_{11}^{(1)} \oplus \beta_{11}^{(2)} \beta_{12}^{(1)} \left( \beta_{21}^{(2)} \beta_{12}^{(1)} \right)^* \beta_{21}^{(2)} \beta_{11}^{(1)} \right] U_{fw}^{(1)} \\ &\oplus \left[ \beta_{11}^{(2)} \beta_{12}^{(1)} \left( \beta_{21}^{(2)} \beta_{12}^{(1)} \right)^* \beta_{22}^{(2)} \oplus \beta_{12}^{(2)} \right] U_{bw}^{(2)}. \end{aligned} \tag{7}$$

Similarly, from the dynamics of system 1, we have
$$Y_{bw}^{(1)} = \beta_{21}^{(1)} U_{fw}^{(1)} \oplus \beta_{22}^{(1)} U_{bw}^{(1)}.$$

Then by replacing $U_{bw}^{(1)}$ by the expression of $Y_{bw}^{(2)}$ in (6), we get
$$Y_{bw}^{(1)} = \beta_{21}^{(1)} U_{fw}^{(1)} \oplus \beta_{22}^{(1)} \left[ \left( \beta_{21}^{(2)} \beta_{12}^{(1)} \right)^* \left( \beta_{21}^{(2)} \beta_{11}^{(1)} U_{fw}^{(1)} \oplus \beta_{22}^{(2)} U_{bw}^{(2)} \right) \right].$$

That is

$$\begin{aligned} Y_{bw}^{(1)} = &\left[ \beta_{21}^{(1)} \oplus \beta_{22}^{(1)} \left( \beta_{21}^{(2)} \beta_{12}^{(1)} \right)^* \beta_{21}^{(2)} \beta_{11}^{(1)} \right] U_{fw}^{(1)} \\ &\oplus \left[ \beta_{22}^{(1)} \left( \beta_{21}^{(2)} \beta_{12}^{(1)} \right)^* \beta_{22}^{(2)} \right] U_{bw}^{(2)} \end{aligned} \tag{8}$$

The result is given by (7) and (8). ∎